\documentclass[11pt]{amsart}

\usepackage[english]{babel}
\usepackage{amssymb}
\usepackage{fourier}
\usepackage{mathrsfs}
\usepackage{enumerate}
\usepackage{color}
\usepackage{ifpdf}
\usepackage{tikz}
\usepackage{tikz-cd}
\usetikzlibrary{decorations.pathmorphing,arrows, intersections}
\usepackage[initials]{amsrefs}
\usepackage{enumitem}
\usepackage{thmtools}
\usepackage{thm-restate}
\usepackage{hyperref}

\usepackage{cleveref}

\usepackage{caption}
\usepackage{subcaption}
\usepackage{comment}

\newcommand{\bbR}{{\mathbf R}}

\newcommand{\calC}{\mathcal{C}}
\newcommand{\calD}{\mathcal{D}}

\newcommand{\PSL}{\operatorname{PSL}}

\newcommand{\Isom}{\operatorname{Isom}}

\newcommand{\setdef}[2]{ \left\{ {#1}\ /\ {#2} \right\} }

\newtheorem{theorem}{Theorem}[section]
\newtheorem{thm}[theorem]{Theorem}
\newtheorem{mthm}{Theorem}

\newtheorem{lemma}[theorem]{Lemma}

\newtheorem{corollary}[theorem]{Corollary}

\newtheorem{proposition}[theorem]{Proposition}

\newtheorem{definition}[theorem]{Definition}

\newtheorem{example}[theorem]{Example}

\newtheorem{remark}[theorem]{Remark}

\numberwithin{equation}{section}

\begin{document}
\title[Marked length spectrum rigidity]{Marked length spectra of Gromov hyperbolic space}
\author{Yanlong Hao}
\address{University of Michigan, Ann Arbor, Michigan, USA}
\email{ylhao@umich.edu}
\keywords{marked length variety, 
marked length pattern,   
 arithmetic surface}
 
 \subjclass[2020]{37D40; 37A20}
\date{\today}
    
\maketitle

\begin{abstract}
    Let $(X,d)$, $(Y, d')$ be two roughly geodesically complete Gromov hyperbolic spaces under comparable isometric actions of $\Gamma$. Assume that the limit set $\Lambda \Gamma=\partial X\partial Y$. If spaces $X$ and $Y$ have the same asymptotic marked length spectrum, meaning that
    \[\lim_{{\ell_{d}([\gamma])\to \infty, [\gamma]\in \calC^*}}\frac{\ell_d(\gamma)}{\ell_{d'}(\gamma)}=1.\] Then $(X,d)$ and $(Y,d')$ are $\Gamma$-equivariantly roughly isometric.
\end{abstract} 
\section{Introduction}
Given a closed Riemannian manifold $(M,g)$ with negative sectional curvatures and its fundamental group $\pi_1(M)=\Gamma$, a unique closed geodesic exists within each free homotopy class that represents a non-trivial conjugacy class of $\Gamma$. Denoting the set of these conjugacy classes of $\Gamma$ as $\mathcal{C},$ we define the \emph{marked length spectrum} map for $g$: 
$$\ell_g: \mathcal{C}\to \mathbb R^+$$
by associating each class in $\mathcal{C}$ with the length of the corresponding closed geodesic. Denote by $\calC^*$ the set of primitive conjugacy classes. Recall that an element $\gamma$ in a group $\Gamma$ is called primitive if $\gamma\neq \eta^n$ for any $\eta\in\Gamma$ and $|n|>1$. 

The \emph{Marked Length Spectrum Rigidity Conjecture} (cf. Burns--Katok \cite{Burns-Katok})
asserts that the Riemannian metric $g$ on $M$ is uniquely determined by the function $\ell_g$ on $\calC_\Gamma$.
More precisely, the conjecture states that if $g_1, g_2$ are two negatively curved Riemannian structures on $M$,
then $\ell_{g_1}=\ell_{g_2}$ only if there is an isometry isotopic to identity $(M,g_1)\cong (M,g_2)$.
This conjecture, motivated by \cite{Guil-Kazh}, was proved in some cases \cite{Otal},
\cite{Croke}, most recently \cite{GL}.
In particular, if $M$ admits a locally symmetric metric, and $g_1$ is such a metric
Hamenst\"{a}dt proved the MLSR conjecture \cite{Ham} (using \cite{BCG}).

A coarsely geometric generalization of the marked length spectrum map is also applicable in this context. Considering a metric space $(X,d)$, along with an isometry  $f: X\to X$, and a fixed base point $o\in X$, we define the  \emph{translation length} of $f$ as 
\[
\ell_d(f)=\lim_{n\to \infty}\frac{d(f^n(o),o)}{n}.
\]
This definition, referring to the translation length, is independent of the choice of the base point $o\in X$. When $(X,d)$, serves as the Riemannian universal cover of $(M,g)$, and $f\in \Gamma=\pi_1(M)$ acts as a deck transformation, the translation length coincides with the marked length spectrum map. Consequently, we refer to the map $\ell_d$ as the marked length spectrum of the isometry group of $X$.

In this context, we have the following theorem:


\begin{mthm}\label{asymptotics1}
Let $(X, d_1)$ and $(Y, d_2)$ be two Gromov hyperbolic spaces under non-elementary isometric actions by a group $\Gamma$. Assume that there is a $\Gamma$-equivalent continuous map from $\Lambda_X \Gamma$ to $\Lambda_Y \Gamma$, where $\Lambda_X \Gamma$ and $\Lambda_Y \Gamma$ are  the limit sets of $\Gamma$ for these two actions, respectively. If
\[
\lim_{\ell_{d_1}([\gamma])\to \infty, [\gamma]\in \calC^*}\frac{\ell_{d_1}([\gamma])}{\ell_{d_2}([\gamma])}=1,
\]
then $\ell_{d_1}=\ell_{d_2}$ on $\Gamma$.
\end{mthm}
Here, $\Lambda_X \Gamma$ and $\Lambda_Y \Gamma$ are the limit set of the two actions, respectively.

By the work of Otal (\cite{Otal}), Croke (\cite{Croke}), and Hamens\"{a}dt (\cite{Ham}), it leads to the following result.
\begin{corollary}
Let $(M,g_1)$ be a closed negatively curved surface or a closed negatively curved locally symmetric space. Denote the fundamental group by $\pi_1(M)=\Gamma$. For any negatively curved Riemannian metric $g_2$ on $M$, if 
\[
\lim_{\ell_{g_1}([\gamma])\to \infty, [\gamma]\in \calC^*}\frac{\ell_{g_1}([\gamma])}{\ell_{g_2}([\gamma])}=1,
\]
then there is an isomorphism isotopic to identity between $(M,g_1)$ and $(M,g_2)$.
\end{corollary}

Along the way of proving Theorem~\ref{asymptotics1}, we obtained additional results on the structure of the marked length spectrum. We start with spectrally rigid sets.
\subsection{Spectrally rigid set}

We first delve into spectrally rigid sets, in exploring properties of the marked length spectrum. Consider a group $\Gamma$ and two isometric actions of $\Gamma$ on metric spaces $(X,d)$ and $(Y, d')$. As in \cite{cantrell2023sparse}, we say a collection of conjugacy class $E$ \emph{spectrally rigid} if, when $\ell_d([x])=\ell_{d'}([x])$ for all $[x]\in E$, it follows that $\ell_d=\ell_{d'}$ for all $[x]\in \mathcal{C}$. By selecting a representative for each class, we may identify $E$ as a subset of $\Gamma$.

Sawyer \cite{sawyer2020partial} has shown that when the marked length spectrum coincides on a subset of $\Gamma$, whose complement exhibits sub-exponential growth,  it uniquely determines the marked length spectrum of a compact surface. Cantrell and Cantrell-Reyes \cites{cantrell2022manhattan, cantrell2023marked, cantrell2023sparse} have extended this notion to a coarse-geometric setting, yielding remarkable results. Notably, Cantrell has constructed arbitrarily sparse spectrally rigid sets for hyperbolic groups.

In this paper, we extend we extend the results discussed earlier to  
certain subgroups of $\Isom(X)$ where  $X$ is a Gromov hyperbolic space. 

Before stating our results, let us establish some notation. Let $X$ be a Gromov hyperbolic space. We call $X$ \emph{nontrivial} if the Busemann boundary $\partial X$ contains at least 3 points. 
Given two non-elementary actions of $\Gamma$ on $(X,d)$ and $(Y, d')$, respectively, the actions of $\Gamma$ on $(X,d)$ and $(Y,d')$ are called \emph{comparable} if there exists a $\Gamma$-equivariant, continuous homeomorphism from $\Lambda_X\Gamma$ to $\Lambda_Y\Gamma$. 

In this paper, all spaces will be 
Gromov hyperbolic spaces.


In this direction, we have the following theorems.
\begin{mthm}\label{mthm: sparse rigid set}
 Let $f:\mathbb R\to \mathbb R^+$ be a function such that $\lim_{T\to\infty}f(T)=\infty$, and $\Gamma$ be a group. Given two non-elementary or quasi-parabolic comparable actions of $\Gamma$ on two non-trivial Gromov hyperbolic spaces $(X, d)$ and $(Y, d')$, for any point $\xi$ in the limit set $\Lambda_X \Gamma$, there exists a subset $E$ such that if $\ell_d(\gamma)=\ell_{d'}(\gamma)$ for all $\gamma\in E$, then $\ell_d=\ell_{d'}$ on $\Gamma$. Furthermore,  we can choose $E$ such that:
 \begin{enumerate}
     \item\label{1} the growth of $E$ with respect to $d$ (and/or $d'$) satisfies:
 \[
 \#\{\gamma\in E|\ell_d(\gamma)\leq T\}\leq f(T).
 \]
 \item\label{2} The only accumulation point of $E$ on $\partial X$ is $\xi$.
 \item\label{4} If furthermore, the action on $X$ is \textbf{proper}, we may choose $E$ which contains only primitive elements.
\end{enumerate}
\end{mthm}
\begin{remark}
    Without proper conditions, condition \ref{4} cannot be achieved. Here is an example: for any subgroup $\Gamma\in \PSL(2,\mathbb R)$, we denote $Root(\Gamma)$ the subgroup of $\PSL(2,\mathbb R)$ generated by all $g$ with $g^2\in \Gamma$. Let $G$ be a Lattice in $\PSL(2,\mathbb R)$, the group $A=\bigcup_{i=0}^\infty Root^i(G)$ is a countable group without any primitive element and acts on the hyperbolic plane with limit set the whole circle.
\end{remark}


While the methodology proposed in \cite{cantrell2023sparse} is expected to apply to cofinite volume actions, certain challenges arise. Firstly, as noted in \cite{cantrell2023marked}*{page 2, last paragraph}, geodesic flow and Busemann cocycle are in general not available outside the Riemannian or low dimensional setting. In this paper, we introduce a rough version of the Busemann cocycle for Gromov hyperbolic space, which serves as a valuable tool in numerous cases. Secondly, the generalization of their work to find spectrally rigid subsets of general subgroups of $\Isom(X)$ remains uncertain. 
Hence, we propose an alternative construction for finding spectrally rigid subsets that also yield an arbitrarily sparse spectrally rigid set.

When applying Theorem~\ref{mthm: sparse rigid set} to a hyperbolic group acts on itself with word metric, we recover the results of \cite{cantrell2023sparse}.
\medskip

On the other hand, our methods are sufficiently flexible that they apply to multiplicative estimates; see Theorem~\ref{mthm: sparse multiplicative rigid set} below, which establishes the applicability of our methods to multiplicative estimates. 

The foundational work of Butt \cites{butt2022quantitative, butt2022approximate} serves as the inspiration for Theorem~\ref{mthm: sparse multiplicative rigid set}. Butt ingeniously proposed that for two distinct Riemannian metrics of negative curvature on a closed manifold of dimension at least 3, a closely matched marked length spectrum implies their volumes are closely aligned. Additionally, she showed that, under certain geometric conditions, the values of the marked length spectrum on a sufficiently large finite set can effectively determine the metric. In this work, we extend Butt's results regarding the marked length spectrum map to encompass all non-elementary or quasi-parabolic actions on Gromov hyperbolic spaces, but for arbitrarily sparse subsets of $\Gamma$.
Given the lack of control over the geometry in our setting, we do not anticipate achieving the same results for a finite subset. 

\begin{mthm}\label{mthm: sparse multiplicative rigid set}
  Let $f:\mathbb R\to \mathbb R^+$ such that $\lim_{T\to\infty}f(T)=\infty$ and $\Gamma$ be a group.  Given two non-elementary or quasi-parabolic comparable actions of $\Gamma$ on two nontrivial Gromov hyperbolic spaces $(X,d)$ and $(Y, d')$, and a point $\xi\in\Lambda_X \Gamma$, there exist a subset $E$ of $\Gamma$ such that: 
  \begin{enumerate}
      \item\label{11} for all $\gamma\in \Gamma$,
  \[\frac{\ell_d(\gamma)}{\ell_{d'}(\gamma)}\geq \liminf_{\ell_d(\eta)\to \infty, \eta\in E}\frac{\ell_d(\eta)}{\ell_{d'}(\eta)},\]  
  \[\frac{\ell_d(\gamma)}{\ell_{d'}(\gamma)}\leq \limsup_{\ell_d(\eta)\to \infty, \eta\in E}\frac{\ell_d(\eta)}{\ell_{d'}(\eta)}.\] 
  \item the growth of $E$ with respect to $d$ (or/and $d'$) satisfying:
 \[
 \#\{\gamma\in E|\ell_d(\gamma)\leq T\}\leq f(T).
 \] 
 \item The only accumulation point of $E$ on $X$ is $\xi$.
 \item Additionally, if the action on $X$ is proper or $\Gamma$ is a torsion-free word hyperbolic group, then E contains only primitive elements.
   \end{enumerate}
\end{mthm}

We term a subset $E$, \emph{detecting the ratio of the marked length spectrum} if it satisfies \ref{11} in Theorem~\ref{mthm: sparse multiplicative rigid set}. 
\subsection{Marked length spectrum rigidity} Now we turn to marked length spectrum rigidity in the coarsely geometric setting.
First, Furman \cite{F} established the following rigidity theorem:
\begin{thm}
    Consider a countable hyperbolic group $\Gamma$, and let $d_1$ and $d_2$ be two roughly geodesic left invariant metrics on $\Gamma$, which are quasi-isometric to a word metric, with a finite set of generators. If $\ell_{d_1}=\ell_{d_2}$ on $\Gamma$, meaning that for all for all $\gamma\in \Gamma$, $\ell_{d_1}(\gamma)=\ell_{d_2}(\gamma)$, then there exist a constant $C$ such that $|d_1-d_2|<C$.
\end{thm}

 Fujiwara \cite{fujiwara2015asymptotically} built upon a strategy initially proposed by Krat \cite{krat2001pairs}, extended this result to toral relatively hyperbolic groups. In 2022, Nguyen and Wang \cite{nguyen2022marked} further confirmed this result for all relative hyperbolic groups. It's important to note that all of these results apply in situations where the metrics are quasi-isometric to a word metric. However, this raises the question: what about the general case? We attempt to answer this question in a simple scenario—the case of isometric actions on Gromov hyperbolic spaces.

Our initial focus is on negatively curved non-compact manifolds with finite volume and bounded geometry. We present the following theorem:
\begin{mthm}\label{mthm: manifold}
    Let $M$ be a manifold of dimension $d\geq 2$. Suppose $g_1$ and $g_2$ are two complete, finite volume, Riemannian metrics on $M$ with sectional curvature satisfying $-a^2\leq \kappa\leq -b^2<0$. If the marked length spectrum of $g_1$ and $g_2$ coincide, then the universal covers of the two metrics are roughly isometric. 
\end{mthm}

\medskip

On the other hand, there's a genuine desire to comprehend the partial rigidity of the marked length spectrum, a concept that has attracted many works in recent years. Partial rigidity refers to the phenomenon where certain subsets of information about a manifold's marked length spectrum can determine or constrain its entire marked length spectrum. For example, Katok \cite{katok1988four} showed that given a closed manifold $M$, for any homology class $[c]\in H_1(M,\mathbb Z)$, two negatively curved metrics in a fixed conformal class possess the same marked length spectrum if and only if the marked length spectrum map is in agreement for all loops representing $[c]$. Gogolev and Rodriguez Hertz \cite{gogolev2020abelian} further advanced our understanding by strengthening Katok's result, showing that it holds even when the condition 'in a fixed conformal class' is dropped. See Remark 6.5 in \cite{gogolev2020abelian}. In our previous work \cite{hao2022marked}, we showed that it is enough to consider the restriction of the marked length spectrum to a subgroup with the same limit set as the fundamental group.

\medskip

One of the principal focuses of this paper is to extend these results to certain groups that act isometrically on a `good' Gromov hyperbolic space.

Before delving into this further, let us establish some notation first. Consider $(X,d)$, a $K$-roughly geodesic, $\delta$-Gromov hyperbolic space. We say that $X$ \emph{$K$-visual} if, for every point $p\in X$ and every point $x\in \partial X$, there exists a $K$-roughly geodesic path starting from $p$ and converging to $x$. We define $(X,d)$ as \emph{$K$-roughly complete} if, for any two points $p$, $q\in X$, there exists an infinite $K$-rough geodesic path passing $p$, $q$ and converging to some $x\in \partial X.$ Note that a proper roughly complete space is always visual.
\begin{definition}
    A metric space $(X,d)$ is called roughly geodesically complete if there exist  $K$, $\delta>0$ such that the following are satisfied:
    \begin{enumerate}
        \item $(X,d)$ is $K$-roughly geodesic, $\delta$-Gromov hyperbolic space;
        \item $(X,d)$ is $K$-visual;
        \item $(X,d)$ is $K$-roughly complete.
    \end{enumerate}
\end{definition}
\begin{example}
    Let $\Gamma$ be a non-elementary hyperbolic group. Any word metric with finitely many generators is a roughly geodesically complete Gromov hyperbolic space.
\end{example}

Given a Gromov hyperbolic space, 
 a non-elementary or quasi-parabolic subgroup $\Gamma$ of $\Isom(X,d)$ is \emph{geometrically dense} if the limit set of $\Gamma$, $\Lambda \Gamma$ is the same as $\partial X$. This notation is inspired by the work of Osin \cite{osin2017invariant}.  

There are many examples of geometrically dense actions. We highlight three of our major concerns.
\begin{example}
    A normal (or commensurated) subgroup of a hyperbolic group acts on the hyperbolic group with word metric.
\end{example}
\begin{example}
    Consider a dense (in Hausdorff topology) finite rank free subgroup $F$ of $\PSL_2(\mathbb R)$. The $F$-action on the hyperbolic plane $H^2$ is geometrically dense. 
\end{example}

Now, we present the main result of the marked length spectrum rigidity.
\begin{mthm}\label{mthm}
Let $\Gamma$ be a group. Given two comparable actions of $\Gamma$ on nontrivial roughly geodesically complete Gromov hyperbolic spaces $(X,d)$ and $(Y,d')$ which are geometrically dense, if $\ell_d=\ell_{d'}:\Gamma\to \mathbb{R}^{\geq 0}$, then there exist a $\Gamma$-equivariant roughly isometry $\phi: X\to Y$, i.e., there exist a constant $C>0$, such that for all $p$, $q\in X$, 
\[|d'(\phi(p),\phi(q))-d(p,q)|<C.\]
\end{mthm}

\begin{remark}
It is evident that Theorem~\ref{mthm: manifold} is a consequence of Theorem~\ref{mthm}. Hence, we will focus our attention on proving the latter. 
\end{remark}

Note that the very recent work of Wan, Xu, and Yang \cite{wan2024marked} shows that for cocompact or cusp-uniform actions on Gromov hyperbolic spaces, the results of Theorem~\ref{mthm} hold. Since in Theorem\ref{mthm}, we only require $\Gamma$ to be geometrically dense, suggesting that our results have broader applicability beyond the settings considered in \cite{wan2024marked}. 

The proof strategy of Theorem~\ref{mthm} is similar to Furman's work \cite{F}. In Furman's work, the key insight is that when two spaces $X$ and $Y$ have the same marked length spectrum, their cross-ratios on the boundary coincide. To leverage this insight, we introduce the concept of a rough cross-ratio on the boundary and demonstrate its close relationship with the marked length spectrum. Unlike Furman's definition, which is a measurable function, our rough cross-ratio is defined everywhere on the boundary, albeit not continuously. 

With the rough cross-ratio established, we proceed to consider the space of metrics $\calD^R$ on the boundary that share the same rough cross-ratio and exhibit bounded geometry. Utilizing a generalization of the barycenter method employed in \cite{biswas2015moebius}, we show that $X$ and $Y$ are roughly isometric to $\calD^R$. Consequently, Theorem~\ref{mthm} follows.

It's worth noting that the space constructed in the second step serves as a hyperbolic extension of the boundary. While other methods, such as comparing hyperbolic spaces of constant curvature, could achieve similar results. Previous work in this direction includes Paulin \cite{paulin1996groupe}, Bonk and Schramm \cite{bonk2011embeddings}, Buyalo and Schroeder \cite{buyalo2007elements}, and Jordi \cite{jordi2010interplay}, among others. we opt for the barycenter method due to its close alignment with the first step of our proof strategy. For a more in-depth discussion on these techniques, refer to section \ref{barycenter}.

\medskip

Note that in Katok \cite{katok1988four} and Gogolev and Rodriguez Hertz \cite{gogolev2020abelian}, the subsets considered are not are not necessarily subgroups. We extend this consideration to cosets of a subgroup. This generalization allows for a broader application of the marked length spectrum rigidity theorem, encompassing a wider range of group actions.
\begin{mthm}\label{Mthm: coset}
 Let $\Gamma$ be a group. Given two comparable actions of $\Gamma$ on roughly geodesically complete Gromov spaces $(X,d)$ and $(Y, d')$ which are geometrically dense, assume that there exists a subgroup $\Gamma'$ and an element $h\in \Gamma$ so that the limit set $\Lambda \Gamma'=\Lambda \Gamma$ for both actions, and $\ell_d=\ell_{d'}$ on the left coset $h\Gamma'$ or the right coset $\Gamma' h$, then $\ell_d=\ell_{d'}$ on $\Gamma$, and there exist a $\Gamma$-equivariant roughly isometry $\phi: X\to Y$, i.e., there exist a constant $C>0$, such that for all $p$, $q\in X$, 
\[|d'(\phi(p),\phi(q))-d(p,q)|<C.\]    
\end{mthm}
\begin{remark}
    Since the two actions are comparable, the condition $\Lambda \Gamma'=\Lambda \Gamma$ for one action is sufficient..
\end{remark}
\begin{remark}
     It is crucial for the space the group acts on to be roughly geodesic because, in general, a $\Gamma$-orbit is not roughly geodesic as a metric space. may not be roughly geodesic as a metric space. For example, consider a non-uniform lattice $\Gamma$ of $\PSL_2(\mathbb R)$ acting on the hyperbolic plane $H^2$. Any $\Gamma$  orbit is not roughly geodesic with the induced metric. However, to prove Theorem~\ref{mthm: sparse rigid set} and Theorem~\ref{mthm: sparse multiplicative rigid set} in the before subsection, it suffices to consider only one $\Gamma$-orbit.
\end{remark}

\medskip

The paper is structured as follows: In section 2, we recall some background and provide some useful lemmas. Section 3 is devoted to the length relation, which is the main idea of this paper. In section 4, we prove the Theorems~\ref{mthm: sparse rigid set} and ~\ref{mthm: sparse multiplicative rigid set}. Then we define and analyze the infinite barycenter map in Section 5.  Section 6 is proof of Theorem~\ref{main}. Finally, we prove  Theorem~\ref{Mthm: coset} in Section 7.

\medskip

The author expresses gratitude to Ralf Spatzier for valuable comments and for drawing attention to the works of Butt and Cantrell. Thanks to Cantrell for many suggestions and comments. Additionally, thanks to Ashot Minasyan for highlighting reference \cite{minasyan2010normal}.
\section{Roughly geodesically complete hyperbolic space}
\subsection{Gromov product and Gromov hyperbolic spaces}
Some good general references for Gromov hyperbolic spaces are \cite{gromov1987hyperbolic} and \cite{bridson2013metric}*{Chapter III.H, Sections 1 and 3}. We recall some basic definitions and results here.

Let $(X,d)$ be a metric space with a base point $o$. The \emph{Gromov product} $<,>_o: X\times X\to \mathbb R^{\geq 0}$ is given by 
\[<p,q>_o=\frac{1}{2}(d(p,o)+d(q,o)-d(p,q)).\]
Let $\delta>0$, $(X,d)$ is called \emph{$\delta$-Gromov hyperbolic} if for all $o$, $p$, $q$, $r\in X$, 
\[
<p,q>_o\geq \min\{<p,r>_o,<r,q>_o\}-\delta.
\]
Note that if there exists some $o\in X$ such that the above inequality holds for all $p$, $q$, $r\in X$, then $X$ is $2\delta$-Gromov hyperbolic. See \cite{gromov1987hyperbolic}*{1.1 B}. We say $(X,d)$ a Gromov hyperbolic space if it is $\delta$-Gromov hyperbolic space for some $\delta$. 

Let $(X,d)$ be a hyperbolic space. A sequence $(x_n)$ in $X$ \emph{converges at infinity} if $<x_i,x_j>_o\to \infty$ as $i$, $j\to \infty$. Two such sequence $(x_n)$ and $(y_n)$ are said to be \emph{equivalent} if $<x_i,y_j>_o\to \infty$ as $i$, $j\to \infty$. The equivalence class of $(x_n)$ is denoted $\lim x_n$, and the set of equivalence classes is denoted $\partial X$, the \emph{Gromov boundary} of $X$.

The Gromov product of a $\delta$-hyperbolic space $X$ can be extended to $\bar X=X\cup \partial X$ by:
\[
<x,y>_o=\sup\liminf_{i,j\to\infty}<x_i,y_j>,
\]
where the supremum is taken over all sequences $(x_i)$ and $(y_j)$ in $X$ such that $x=\lim x_i$ and $y=\lim y_j$. 

\begin{remark}
    The supremum and $\liminf$ lead to difficulty in computing. However, it is sufficient to consider any sequence converging to $x$ and $y$ when only a rough estimation is needed. We have, for all $x$, $y\in \partial X$, and all sequence $(x_i')$  and $(y_j')$ in $X$ with $x=\lim x_i'$ and $y=\lim y_j'$, 
\begin{equation}\label{sequence free}
    <x,y>_o-2\delta\leq \liminf_{i,j}<x_i',y_j'>_o\leq <x,y>_o.
\end{equation}
See \cite{bridson2013metric}*{pp 432-433}.
\end{remark}
Note that the extended Gromov product is not continuous in general. Nevertheless, is is a roughly continuous. Recall (see \cite{bridson2013metric}*{pp 433}) that for all $x$, $y$, $z\in\bar X$, we have 
\begin{equation}\label{hyperolic inequality}
 <x,y>_o\geq \min\{<x,z>_o, <y,z>_o\}-2\delta,   
\end{equation}
and the topology on $\bar{X}$ is given by $x_n\to x$ if $<x_n,x>_o=\infty$ or $x_n$ converge to $y$ in $X$, it follows for any sequence $(y_n)$ in $\bar{X}$ with $y_n\to y$ and $z\in \bar{X}$, we have 
\begin{equation}\label{rough continuous}
    |\liminf<y_n,z>_o-<y,z>_0|\leq 2\delta.
\end{equation}
When $y\in X$, this follows by taking the limit of the triangle inequality.

Finally, we define the (logarithmic) \emph{rough cross-ratio}. Let $\partial^{(4)}X$ be the set of pairwise different 4-tuples in $\partial X$. The rough cross-ratio is defined by the map 
\begin{equation}\label{rough cross-ratio}
(x,y;z,w)=<x,y>_o+<z,w>_o-<z,y>_o-<x,w>_o.
\end{equation}

Note that when $(X,d)$ is a CAT(-1) space, our definition is proportional to the logarithm of the standard cross-ratio defined by Otal \cite{otal1992geometrie}.  
\subsection{Quasi-isometries, roughly geodesic space} We follow the definition in \cite{bonk2011embeddings}. A map $f:(X,d)\to (Y,d')$ of metric spaces is called a $(K,C)$-\emph{quasi-isometric embedding} if 
\[
\frac{1}{K}d'(f(p),f(q))-C\leq d(p,q)\leq Kd'(f(p),f(q))+C,
\] 
for all $p$, $q\in X$. A $(K, C)$-quasi-isometric embedding is called a \emph{quasi-isometry}, if the image is \emph{cobounded}, i.e., there exist constant $L$ such that for all point $r\in Y$, there exist $p\in X$, with $d'(r,f(p))<L$.

$(1, C)$-quasi-isometric embeddings and $(1,C)$-quasi-isometries are called $C$-\emph{roughly isometric embeddings} and $C$-\emph{rough isometries}, respectively. A $C$-\emph{rough geodesic} is a $C$-roughly isometric embedding of an interval. A metric space is called $C$-roughly geodesic if any two points can be connected by a $C$-rough geodesic. 

Let $\tau(t)$ be a rough geodesic, for any $t_1<t_2<t_3$, we have
\begin{equation}\label{Gromov 1; rough geodesic}<\tau(t_1),\tau(t_3)>_{\tau(t_2)}\leq \frac{3}{2}C,
\end{equation}
\begin{equation}\label{Gromov 2: rough geodesic}
<\tau(t_2),\tau(t_3)>_{\tau(t_1)}\geq |t_2-t_1|-\frac{3}{2}C.
\end{equation}

A rough geodesic space is called \emph{roughly complete} if every rough geodesic extends to a bi-infinite rough geodesic.

Let $\Gamma$ be a group acting isometrically on a hyperbolic space $X$, the limit set of $\Gamma$, $\Lambda_X\Gamma$ is the set $\overline{\Gamma\cdot o}\cap \partial X$ for any point $o\in X$.







    

\subsection{Busemann cocycle (rough Busemann cocycle)}\label{Busemann cocycle}
Since its initial definition in \cite{busemann2012geometry}, the Busemann function has found applications in various problem domains. Notably, several applications involve its related counterpart, the Busemann cocycle, as discussed in \cites{benoist2016central, nica2013proper, nica2019two}. We first recall the definition of the Busemann function. Let $(X,d)$ be a Gromov-hyperbolic space.  
For points $o,p\in X$ and $x\in \partial X$, the function 
$B:X\times X\times \partial X\rightarrow \bbR$ is defined by:
 \begin{equation}\label{eq: Busemann function}
    B_{o,p}(x)=\sup\liminf_{q_i\rightarrow x}d(o, q_i)-d(p,q_i)=2<p,x>_o-d(o,p),
\end{equation}
where $(q_i)$ runs all the sequence converge to $x$.

\begin{definition}
    Let $(X,o,d)$ be a pointed Gromov hyperbolic space. The (rough) Busemann cocycle is the function:
    \[c:\Isom(X)\times \partial X\to \mathbb R,\]
    \[
    c(\gamma,x)=B_{o,\gamma o}(x).
    \]
    The extended Busemann cocycle is a similar formulation:
    \[c:\Isom(X)\times X\to \mathbb R,\]
    \[
    c(\gamma,p)=d(o,p)-d(\gamma o,p)=2<\gamma o,p>_o-d(o,\gamma o).
    \]
\end{definition}
\begin{remark}
    In earlier works on Gromov hyperbolic spaces, the Busemann cocycle is defined as the Radon-Nikodym derivative of the Patterson-Sullivan measure, constituting a cocycle. However, in this context, they are merely measurable functions, making it challenging to compare different Busemann cocycles without information about their Patterson-Sullivan measures. The definition adopted here maintains close to the Radon-Nikodym derivative, resulting in a rough cocycle that shares numerous properties. Importantly, the new functions are defined everywhere. 
\end{remark}

Here, we enumerate some properties of Busemann cocycles, with concise proofs provided for the reader's convenience. Most of these properties are well-established, particularly when the Busemann function is continuous  \cite{nica2019two}. 
\begin{proposition}\label{almost continuously}
    The Busemann cocycle $c$ is almost continuous. For any sequence $(y_i)$ in $\bar{X}$ converging to $x\in \partial X$,  and for any $\gamma\in\Isom(X)$,
    \[
    |\liminf_{i\to\infty}c(\gamma,y_i)-c(\gamma,x)|\leq 4\delta.
    \]
\end{proposition}
\begin{proof}
    This follows directly from the definition of the Busemann function and equation~(\ref{rough continuous}).
\end{proof}
\begin{proposition}\label{almost cocycle}
    For any $\gamma$, $\gamma'\in \Isom(X)$ and all $x\in \partial X$, we have
    \[
    |c(\gamma\gamma',x)-c(\gamma,\gamma'x)-c(\gamma',x)|\leq8\delta.
    \]
\end{proposition}
\begin{proof}
    Choose $\epsilon>0$. Select a sequence $(q_i)$ in $X$ converging to $x$ such that \[\lim_{i\to \infty}c(\gamma\gamma',q_i)-\epsilon\leq c(\gamma\gamma',x)\leq\lim_{i\to \infty}c(\gamma\gamma',q_i)+\epsilon.\] Note that $c(\gamma\gamma',q_i)=c(\gamma,\gamma'q_i)+c(\gamma',q_i)$ for all $i$. By proposition \ref{almost continuously}], 
    \[
    |\lim_{i\to \infty}c(\gamma,\gamma'q_i)-c(\gamma, \gamma'x)|\leq 4\delta,
    \]
     \[
    |\lim_{i\to \infty}c(\gamma',q_i)-c(\gamma', x)|\leq 4\delta.
    \]
   Since $\epsilon$ is arbitrary, the result follows.
\end{proof}
A similar argument (setting proper sequence converge to $x$ and $y$) leads to:
\begin{proposition}\label{B-cocycle}
    For all $\gamma\in \Isom(X)$ and $x$, $y\in \partial X$,
    \[
    |<x,y>-<\gamma x,\gamma y>-c(\gamma, x)-c(\gamma,y)|\leq 4\delta.
    \]
\end{proposition}
\begin{proposition}\label{Prop: stable length}
   For all hyperbolic elements $\gamma\in\Isom(X)$, $c(\gamma,\gamma^\pm)=\pm\ell_d(\gamma)$. 
\end{proposition}
\begin{proof}
    First, consider $c(\gamma,\gamma^+)$. By the fact $\gamma\gamma^+=\gamma^+$, and Proposition~\ref{almost cocycle}, we have:
    \begin{equation}\label{eq: almost cocycle}
      |c(\gamma^n,\gamma^+)-nc(\gamma,\gamma^+)|\leq 8(n-1)\delta.  
    \end{equation}
    
Now, take any sequence $(q_i)$ in $X$ converging to $\gamma^+$. Then, by proposition~\ref{almost continuously},
\[
|\liminf_{i\to \infty}c(\gamma^n,q_i)-c(\gamma^n,\gamma^+)|\leq 4\delta.
\]
Notice that
\[
\begin{array}{llr}
c(\gamma^n,q_i)  &=d(o,q_i)-d(\gamma^no,q_i)& \text{(Definition)} \\
     &=d(o,\gamma^no)-2<o,q_i>_{\gamma^no}&\\
     &=d(o,\gamma^no)-2<\gamma^{-n}o,\gamma^{-n}q_i>_o.&(\gamma\in\Isom(X))
\end{array}
\]
 Hence by equation~\eqref{rough continuous},
 \[
 |\liminf_{i\to \infty}c(\gamma^n,q_i)-d(o,\gamma^n o)-2<\gamma^{-n}o,\gamma^+>_0|\leq 4\delta.
 \]
 Since $\gamma^{-n}o$ converge to $\gamma^-\neq \gamma^+$, we see that $<\gamma^{-n}o,\gamma^+>_o$ is bounded. Therefore,
 \[\lim_{n\to \infty}\frac{c(\gamma^n,\gamma^+)}{n}=\lim_{n\to \infty}\frac{d(o,\gamma^n o)}{n}=\ell_d(\gamma).\]
 It follows from equation~(\ref{eq: almost cocycle}), $|c(\gamma,\gamma^+)-\ell_d(\gamma)|\leq 8\delta$. 

 For $c(\gamma,\gamma^-)$, notice that by Proposition~\ref{B-cocycle}, 
 \[
 |c(\gamma,\gamma^+)+c(\gamma,\gamma^-)|\leq 4\delta.
 \]
We have, $c(\gamma,\gamma^-)+\ell_d(\gamma)\leq 12\delta$.

Applying this equation to $\gamma^n$, we have
\[
|c(\gamma^n,\gamma^\pm)\mp\ell_d(\gamma^n)|=n|c(\gamma,\gamma^\pm)\mp\ell_d(\gamma)|\leq 12\delta.
\]
Therefore, $c(\gamma,\gamma^\pm)=\pm\ell_d(\gamma)$. This completes the proof.
 
\end{proof}
\subsection{Barycenter method and relative Busemann function}\label{barycenter}
Following Bourdon \cite{bourdon1996birapport},  Biswas \cite{biswas2015moebius} showed that any Moebius homeomorphism
$f : \partial X \to \partial Y$, between boundaries of proper, geodesically complete CAT(-1) spaces
$X$, $Y$, extends to a $(1, \log 2)$-quasi-isometry $F : X \to Y$. Biswas \ further generalized this result \cites{biswas2017hyperbolic, biswas2021quasi} in a more comprehensive framework. We outline some key ideas here; for a detailed exposition, refer to \cite{biswas2021quasi}.

Consider a compact metrizable space $Z$, where an antipodal function is defined as a continuous function $\phi: Z \times Z \to [0, 1]$, satisfying symmetry, positivity, and the antipodal property such that for every $\xi\in Z$,  there exists an antipode $\eta$ with $\phi(\xi, \eta) = 1$. An antipodal space is a compact Hausdorff space $Z$ with an antipodal function $\phi$. Notably, the boundary of a proper, geodesically complete, boundary-continuous Gromov hyperbolic space equipped with a visual quasi-metric serves as an example of an antipodal space. 

Given an antipodal space $(Z,\phi)$, define the corresponding \emph{Moebius space} $M(Z)$ to be the set of Moebius antipodal functions on $Z$, i.e., antipodal functions share the same cross-ratio as the that of the given antipodal function on $Z$.
A metric is defined on $M(Z)$ by 
\[d_M(\phi_1,\phi_2)=\left\lVert\log \frac{d\phi_1}{d\phi_2}\right\rVert_\infty,\]
where $\frac{d\phi_1}{d\phi_2}$ is the derivative of $\phi_1$ with respect to $\phi_2$, and $\left\lVert\cdot\right\rVert$ denotes the supremum norm, see \cite{biswas2021quasi}*{Lemma 2.2}. 

One notable result from \cite{biswas2021quasi} asserts that for a proper, geodesically complete, boundary-continuous Gromov hyperbolic space $X$, $X$ can be isometrically embedded into the space  $(M(\partial X), d_M)$ as a cobounded subset.

We call $\log \frac{d\phi_1}{d\phi_2}$ the \emph{relative Busemann function}, as it represents the difference in the Busemann function when both metrics are negatively curved Riemannian. This concept is then generalized to a rough geodesic setting and agrees, up to bounded error, with the previous definition.

We will have a similar construction in Section~\ref{Rough corss-ratio}, based on the relative Busemann function.
\subsection{Some lemmas}
    Here we recall some results. All of them are probably well-known to experts. However, we did not find precise statements for general actions.
\begin{lemma}\label{dense orbit}
 Let $\Gamma$ be a group. Consider a non-elementary isometric action of $\Gamma$ on a non-trivial $\delta$-Gromov hyperbolic space $(X,d)$. Then, the set of pairs of fixed points of hyperbolic elements is dense in the complement of
the diagonal $(\Lambda_X \Gamma \times \Lambda_X \Gamma)\setminus \Delta$.
\end{lemma}
The proof of Lemma 5.5 in \cite{hamenstadt2009rank} applies here. Note that, only the dynamics of the action on the boundary are used in the proof there. 

\begin{lemma}\label{discete orbit}
 Let $\Gamma$ be a group. Consider a non-elementary isometric \textbf{proper} action of $\Gamma$ on a non-trivial $\delta$-Gromov hyperbolic space $(X,d)$. Then, for any $T>0$, the set of pairs of fixed points of hyperbolic elements with transitive length less than $T$ is discrete in the complement of
the diagonal $(\Lambda_X \Gamma \times \Lambda_X \Gamma)\setminus \Delta$.
\end{lemma}
\begin{proof}
    Assume there is a $T$, and a pair of distinct points $(\xi, \eta)$, which is a limit points of a sequence of pair of fixed point $(\gamma_n^+,\gamma_n^+)$ with $\ell_d(\gamma_n)<T$.  Define $M=<\xi,\eta>_o$. Since the Gromov product is almost continuous, Fro sufficient large $n$, we have 
    \[|<\gamma_n^+,\gamma_n^->-M|\leq2\delta +1.\]
    By definition, 
    \[<\gamma_n o,\gamma_n^->_0=\liminf_{x_j\to\gamma_n^-}\frac{1}{2}(d(o,\gamma_n o)+d(o,x_j)-d(\gamma_n o,x_j))=\frac{1}{2}(d(o,\gamma o))+c(\gamma_n,\gamma_n^-)).\]
    Noting that $c(\gamma_n,\gamma_n^-)=-\ell_d(\gamma_n)$ and similarly, $<\gamma_n o,\gamma_n^+>_o=\frac{1}{2}(d(\gamma_n o, o)+\ell_d(\gamma_n))$, by the hyperbolic inequality \ref{hyperolic inequality}, we have
    \[
    <\gamma_n^+,\gamma_n^->_o\geq \min\{<\gamma_n^+,\gamma_n o>_o, <\gamma_n^-,\gamma_n o>_o\}-2\delta.
    \]
    Thus, for sufficiently large $n$, 
    $\frac{1}{2}(d(\gamma_n o,o)-T)-2\delta\leq M+2\delta +1.$
     Consequently, 
 $d(\gamma_n o,o)\leq 2M+6\delta+T+2$. 

    Since the action is proper, there are only finitely many choices of such $\gamma_n$. Thus the sequence stabilizes, and the results follow.
\end{proof}
\section{Main idea: length relation}
In this section, we present the central method employed in this paper, which focuses on the relation of stable length and the Busemann function. Before delving into the details, we find it beneficial to recall a pertinent lemma. Although it is not exactly a geometric boundary in this context, the second part of the proof of the lemma still applies.

\begin{lemma}\label{S}\cite{hao2022marked}*{Lemma 3.8}
    Let $\Gamma$ be a group. Consider a non-elementary isometric action of $\Gamma$ on a non-trivial $\delta$-Gromov hyperbolic space $(X,d)$. Let $\gamma$, $\eta\in \Gamma$ with $\gamma$ acts as a hyperbolic element on $X$. If $\eta\gamma^+\neq\gamma^-$, then, for $n$ big enough, $\eta\gamma^n$ act
hyperbolic on $X$, and we have:
\[\lim_{n\to \infty}(\eta\gamma^n)^+=\eta\gamma^+,\]
\[\lim_{n\to -\infty}(\eta\gamma^n)^-=\gamma^-.\]
\end{lemma}
With this lemma in mind, we proceed to outline the primary relation mentioned at the start of the section. 
\begin{proposition}\label{main ralation}
    Let $\Gamma$ be a group. Given a non-elementary isometric action of $\Gamma$ on a non-trivial $\delta$-Gromov hyperbolic space $(X,d)$, let $\gamma\in \Gamma$ act as a hyperbolic isometry on $X$, and let $c$ be the rough Busemann cocycle for the action. For any $\eta\in \Gamma$ such that $\eta\gamma^+\neq \gamma^-$, we have, 
    \[
    |c(\eta,\gamma^+)+n\ell_d(\gamma)-\ell_d(\eta\gamma^n)|\leq |<\gamma^+,(\eta\gamma^n)^->-<\eta\gamma^+, (\eta\gamma^n)^->|
+4\delta.    \]
    In particular,
    \[
    \limsup_{n\to \infty}|c(\eta,\gamma^+)+n\ell_d(\gamma)-\ell_d(\eta\gamma^n)|\leq |<\gamma^+,\gamma^->_d-<\eta\gamma^+,\gamma^->_d|+12\delta
    \]
    
\end{proposition}
\begin{proof}
    From section~\ref{Busemann cocycle}, we have the following properties:
    \begin{enumerate}
        \item $c(\gamma,\gamma^{\pm})=\pm\ell_d(\gamma).$ (Proposition~\ref{Prop: stable length})
        \item $|< x, y>_d-<\zeta x,\zeta y>_d-c(\zeta,x)-c(\zeta,y)|\leq 4\delta.$ (Proposition~\ref{B-cocycle})
    \end{enumerate}
    for all $x,y\in \partial X$, $\zeta\in \Gamma$, and $\gamma\in \Gamma$ which is a hyperbolic element. 

    Choose any $\gamma\in \Gamma$ a hyperbolic elements on $X$ with fixed points $\gamma^+$ and $\gamma^{-}$. Let $\eta\in \Gamma$ be any element such that $\eta \gamma^+\neq \gamma^-$. By Proposition~\ref{almost cocycle}, we have:
\begin{equation}
|c(\eta \gamma^n,\gamma^+)-c(\eta,\gamma^n\gamma^+)-c(\gamma^n,\gamma^+)|\leq 4\delta
\end{equation}
for all $n$. Hence, in view of Proposition~\ref{Prop: stable length},
\begin{equation}\label{eq: 1}
|c(\eta \gamma^n,\gamma^+)-c(\eta,\gamma^+)-n\ell_d(\gamma)|\leq 4\delta.
\end{equation}

Now, for sufficiently large $n$, $\eta\gamma^n$ are hyperbolic according to Lemma~\ref{S}. Hence, we have:
\begin{equation}\label{eq: 2}
    c(\eta\gamma^n, \gamma^+)-\ell_d(\eta\gamma^n)=c(\eta \gamma^n, \gamma^+)+c(\eta \gamma^n, (\eta\gamma^n)^-.
\end{equation}

By Proposition~\ref{B-cocycle}, we have:
\begin{equation}\label{eq: 3}
|<\gamma^+,(\eta\gamma^n)^->_d-<\eta \gamma^n \gamma^+, \eta \gamma^n (\eta\gamma^n)^->_d-c(\eta \gamma^n, \gamma^+)-c(\eta \gamma^n, (\eta\gamma^n)^-)|\leq 4\delta.
\end{equation}
Note that $<\eta \gamma^n \gamma^+, \eta \gamma^n (\eta\gamma^n)^->_d=<\eta\gamma^+, (\eta\gamma^n)^->_d$.
Taking the limit as $n\rightarrow +\infty$, using Lemma \ref{S}, equations \eqref{eq: 2} and \eqref{eq: 3} and \eqref{rough continuous}, we obtain:
\begin{equation}\label{eq: 4}
   \limsup_{n\to \infty} |c(\eta,\gamma^+)+n\ell_d(\gamma)-\ell_d(\eta\gamma^n)-<\gamma^+,\gamma^->_d+<\eta\gamma^+,\gamma^->_d|\leq 12\delta.
\end{equation}
The results follow from the triangle inequality.
\end{proof}

Now, we present the first application of the established relation. 

Let $\Gamma$ be a group, and $(X,d)$ a nontrivial $\delta$-Gromov hyperbolic space with a $\Gamma$-action by isometries, which is geometrically dense. According to \cite{gromov1987hyperbolic}*{8.2.D and 8.2.E}, $\Gamma$ is either quasi-parabolic or non-elementary. Let $P(N)$ be the set of all increasing sequences of positive integers, and $\phi:\Gamma\to P(N)$ be any map. Denote the image of $\xi$ by $\phi_\xi$ for all $\xi\in\Gamma$.

\textbf{Case 1: } If $\Gamma$ is quasi-parabolic, and let $\xi\in \Lambda_X \Gamma$. 

\textbf{Subcase 1:} if $\xi$ is the global fixed point. Let $\Gamma$ be a hyperbolic element in $\Gamma$ with $\gamma^+=\xi$. Taking  
\[E_\phi=\{\eta\gamma^{\phi(\eta)(n)}|\eta\in \Gamma, n\in\mathbb N^{+}\}.\]

\textbf{Subcase 2:} if $\xi$ is not the global fixed point. Let $\xi'$ be the global fixed point. Taking disjoint open neighbourhood of $U$ and $V$ of $\xi$ and $\xi'$, respectively, so that there exist $M>0$ such that for all $\xi_1\in U$, $\xi_2\in V$, $<\xi_1,\xi_2><M$. The existence of $U$ and $V$ is guaranteed by the almost continuity of the Gromov product. 

For any hyperbolic element $\eta$, up to conjugation and replace $\eta$ by $\eta^{-1}$, we may assume $\eta^{-1}\in U$. For any $n$, taking a small neighborhood of $\xi$, $W_n$ so that $W_{n+1}\subset W_n$ is a local basis for the point $\xi$. For example, $W_n$ could be an open subset of the set of points $\{\xi_3|<\xi_3,\xi>\geq n\}$. Let $\gamma_{n,\eta}$ be a hyperbolic element so that $(\gamma^+,\gamma^-)\in \eta^{-1}W_n\times V$. In fact, $\gamma^-=\eta$, and such an element can be a conjugation of any hyperbolic element. Let  
\[E_\phi=\{\eta\gamma_{n,\eta}^{\phi(\eta)(n)},\gamma_{n,\eta}|\eta\in \Gamma \text{ is hyperbolic}, n\in\mathbb N^{+}\}.\]

\textbf{Case 2:} If $\Gamma$ is non elementary, and fix $\xi\in \Lambda_X \Gamma$. Choose two more different points $\xi_1$ and $\xi_2\in\Lambda_X \Gamma$ and their open neighborhood $U$, $V\subset V'$ and $W$ of $\xi$, $\xi_1$, and $\xi_2$, respectively, so that the Gromov products of the element in following set has an upper bound $M$: $U\times V$, $U\times W$, $\Lambda_X\Gamma\setminus V'\times V$ and $V'\times W$. Also for any $n$, taking a small neighborhood of $\xi$, $W_n$ so that $W_{n+1}\subset W_n$ is a local basis for the point $\xi$.

Let $\eta\in \Gamma$. If $\eta^{-1}\xi\notin V'$, let $\gamma_{n,\eta}$ be a hyperbolic element so that $(\gamma^+,\gamma^-)\in \eta^{-1}W_n\times V$.This is possible by Lemma~\ref{dense orbit}. If $\eta^{-1}\xi\in V'$. Let $\gamma_{n,\eta}$ be a hyperbolic element so that $(\gamma^+,\gamma^-)\in \eta^{-1}W_n\times W$. Let 
\[E_\phi=\{\eta\gamma_{n,\eta}^{\phi(\eta)(n)}, \gamma_{n,\eta}|\eta\in \Gamma, n\in \mathbb N^+\}.\]

We have the following:
\begin{lemma}\label{thm: marked length rigidity}
 The set $E_\phi$ is spectrally rigid.   
\end{lemma}
\begin{proof}
Let $(Y, d')$ be any nontrivial $\delta$-Gromov hyperbolic space with a comparable $\Gamma$ action. Assume the marked length spectrum agrees on $E_\phi$ for these two actions. We identify $\partial X$ and $\partial Y$ by the homeomorphism $f$. Denote the Busemann cocycle for $X$ and $Y$ by $c_1$ and $c_2$.

\textbf{Case 1: $\Gamma$ is quasi-parabolic.} 

\textbf{Subcase 1: $\xi$ is the global fixed point.} For any $\eta\in \Gamma$, applying Proposition~\ref{main ralation} for $\eta$ and $\gamma$, we have:
\[
\limsup_{n\to\infty}|c_1(\eta,\gamma^+)+\phi_\eta(n)\ell_d(\gamma)-\ell_d(\eta\gamma^{\phi_\eta(n)}|\\
\leq 14\delta,
\]
and 
\[
\limsup_{n\to\infty} |c_2(\eta,\gamma^+)+\phi_\eta(n)\ell_{d'}(\gamma)-\ell_{d'}(\eta\gamma^{\phi_\eta(n)}|\leq 14\delta.
\]
Taking the difference, and note that $\ell_{d}(\eta\gamma^{\phi_\eta(n)})=\ell_{d'}(\eta\gamma^{\phi_\eta(n)})$ by assumption,
we have:
\[
\limsup|c_1(\eta,\gamma^+)-c_2(\eta,\gamma^+)+\phi_\eta(n)(\ell_d(\gamma)-\ell_{d'}(\gamma))|\leq 28\delta.
\]
Hence $\ell_d(\gamma)=\ell_{d'}(\gamma)$ and $c_1(\eta,\gamma^+)-c_2(\eta,\gamma^+)$ are bounded. Since $\gamma^+$ is a fixed point of $\eta$, we have for any $n\in \mathbb Z$, $c_1(\eta^n,\gamma^+)-c_2(\eta^n,\gamma^+)=n(\ell_d(\eta)-\ell_{d'}(\eta))$ is bounded. Therefore, $\ell_d(\eta)=\ell_{d'}(\eta).$

\textbf{Subcase 2: $\xi$ is not the global fixed point.} Let $\xi'$ be the $\Gamma$-fixed point. For any $\eta\in \Gamma$. Since there is a $\Gamma$-equivariant homeomorphism between the two limit sets, if $\eta$ is an elliptic or parabolic isometry of $X$, it is an elliptic or parabolic isometry of $Y$. And $\ell_d(\eta)=\ell_{d'}(\eta)$. 

If $\eta$ is a hyperbolic element for the action on $X$, the same is true for $Y$. 
Applying Proposition~\ref{main ralation} for $\eta$ and $\gamma$, we have that for $n$ big enough:

$|c_1(\eta,\gamma_{n,\eta}^+)+\phi(\eta)(n)\ell_d(\gamma_{n,\eta})-\ell_d(\eta\gamma_{n,\eta}^{\phi(\eta(n))})|$ is less than $|<\gamma_{n,\eta}^+,(\eta\gamma_{n,\eta}^{\phi(\eta)(n)})^->_X-<\eta\gamma_{n,\eta}^+,(\eta\gamma_{n,\eta}^{\phi(\eta)(n)})^->_X|+4\delta+1.$ Same is true for the action on $Y$.

Note that by construction: firstly, $\gamma^+_{n,\eta^m}$ converge to $\eta^{-1}$ uniformly as $m\to \infty$; secondly, $\eta^m\gamma_{n,\eta^m}^+\in W_n$; and finally, $(\eta\gamma_{n,\eta}^{\phi(\eta)(n)})^-$ converge to $\xi'$ (it is a constant sequence, but we do not need this fact). 

We have, for $m$ big enough, and $n$ big enough depending on $m$, all the Gromov products are uniformly bounded for both $X$ and $Y$. By assumption on the marked length spectrum, we have that $|c_1(\eta^m, \gamma_{n,\eta^m}^+)-c_2(\eta^m, \gamma_{n,\eta^m}^+)|$ is uniformly bounded.

By the cocycle identity, $|c_1((\eta^m, \gamma_{n,\eta^m}^+))+c_1(\eta^{-m}, \eta^m \gamma_{n,\eta^m}^+)|\leq 8\delta$. For fixed $m$, choose $n$ big enough, such that: $|c_i(\eta^{-m},\xi)-c_i(\eta^{-m},\xi_1)|\leq 2\delta+1$ for all $\xi_1\in W_n$ and $i=1,2$. We get
$c_1(\eta^{-m},\xi)-c_2(\eta^{-m},\xi)$ are bounded with a bound independent of $m$. Multiplying it by $\frac{1}{m}$, and taking the limit as $m\to \infty$, we have $\ell_d(\eta)=\ell_{d'}(\eta)$.

\textbf{Case 2: $\Gamma$ is non-elementary.} It is very similar to subcase 2 in case 1. We omit it.
\end{proof}
\section{Proof of Theorem~\ref{mthm: sparse rigid set} and \ref{mthm: sparse multiplicative rigid set}}
\subsection{Proof of Theorem~\ref{mthm: sparse rigid set}}
The proof of Theorem~\ref{mthm: sparse rigid set} is an application of Lemma~\ref{thm: marked length rigidity}. First, let us restate Theorem~\ref{mthm: sparse rigid set}.
\begin{theorem}
 Let $f:\mathbb R\to \mathbb R^+$ such that $\lim_{T\to\infty}f(T)=\infty$ and $\Gamma$ be a group. There exists a subset $E\subset \Gamma$ so that: given two non-elementary or quasi-parabolic comparable actions of $\Gamma$ on two non-trivial Gromov hyperbolic spaces $(X,d)$ and $(Y, d')$, if $\ell_d(\gamma)=\ell_{d'}(\gamma)$ for all $\gamma\in E$, then $\ell_d=\ell_{d'}$ on $\Gamma$. Furthermore, for any fixed point $\xi\in\Lambda_X \Gamma$, we can choose $E$ so that:
 \begin{enumerate}
     \item the growth of $E$ with respect to $d$ (or/and $d'$) satisfying:
 \[
 \#\{\gamma\in E|\ell_d(\gamma)\leq T\}\leq f(T).
 \]
 \item The limit set of $E$ contains only $\xi$.
 \item If furthermore, if the action on $X$ is \textbf{proper}, we may choose $E$ which contains only primitive elements.
 \end{enumerate}
\end{theorem}
\begin{proof}
Fix a base point $o\in X$.

 \textbf{Case 1: $\Gamma$ is quasi-parabolic}
 
 \textbf{Subcase 1: $\xi$ is the global fixed point.}
 Let $\Gamma$ be a fixed hyperbolic element such that $\gamma^+$ is a $\Gamma$-fixed point. Note that by Lemma~\ref{S}, for any $\eta\in \Gamma$, for $n$ big enough, $\eta\gamma^{n}$ is hyperbolic, and $|\ell_d(\eta\gamma^{n}-n\ell_{\gamma}|\leq K_{\eta,\gamma}$ by Proposition~\ref{main ralation}. Note that $\eta\gamma^+=\eta\xi=\xi$. Enumerate the element of $\Gamma$ by $\eta_i$, $i\in \mathbb N^+$. For each $i\in \mathbb N^+$,
 choose a sequence $n^i$ so that 
 \begin{enumerate}
     \item $\#\{m|\ell_d(\eta_i\gamma^{n^i(m)})\leq T\}\leq \frac{f(T)}{2^i}$.
     \item $<\eta_i\gamma^{n^i(m)} o,\xi>_d\geq 1000 i^2(\delta+1)$ for all $m>1$.
 \end{enumerate}
 Let $\phi$ be the map $\phi(\eta_i)=n^i$, and set $E=E_\phi.$ Then the only limit point of $E$ is $\xi$, and $E$ is spectrally rigid by Lemma~\ref{thm: marked length rigidity}. By the construction, the growth rate of $E$ is less than $f$.

 \textbf{Subcase 2: $\xi$ is not the global fixed point.} Enumerate the element of $\Gamma$ by $\eta_i$, $i\in \mathbb N^+$. First, choose $\gamma_{n,\eta}$ so that 
 \begin{enumerate}
    \item $\#\{m|\ell_d(\gamma_{n,\eta})\leq T\}\leq \frac{f(T)}{2}$.
     \item $\ell_d(\gamma_{n,\eta})\to \infty$ as $n\to \infty$ for all $\eta\in \Gamma$ acts hyperbolicly.
     \item $<\eta\gamma_{n,\eta}^+,\xi>\geq 10000i^2(\delta+1)$.
 \end{enumerate}
Then notice that for fixed $\eta$, $\ell_d(\eta\gamma_{n,\eta}^m)-m\ell_d(\gamma_{n,\eta})$ is bounded. For each $i\in \mathbb{N}^+$, choose a sequence $n^i$ so that 
  \begin{enumerate}
     \item $\#\{m|\ell_d(\eta\gamma_{m,\eta}^{n^i(m)})\leq T\}\leq \frac{f(T)}{2^{i+1}}$.
     \item $<\eta\gamma_{m.\eta}^{n^i(m)} o,\eta\gamma_{m,\eta}^+)>_d\geq 1000 i^2(\delta+1)$ for all $m>1$.
 \end{enumerate}
 Let $\phi$ be the map $\phi(\eta_i)=m^i$, and set $E=E_\phi$. Firstly, the only limit point of $E$ is $\xi$; secondly, $E$ is spectrally rigid by Lemma~\ref{thm: marked length rigidity}; and thirdly, the growth of $E$ is less than $f$. 
 
\textbf{Case 2: $\Gamma$ is non-elementary.} It is similar to subcase 2 in case 1. We omit it.

Finally, if the action of $\Gamma$ on $X$ is proper, by Lemma~\ref{discete orbit}, every hyperbolic element is a power of a primitive element. Hence in the construction, we may replace all elements with the corresponding primitive elements. Since the $\gamma_{n,\eta}$ can be any element with a pair of fixed points in a prescribed open set. By lemma~\ref{dense orbit}, there are infinitely many choices. And because of Lemma~\ref{discete orbit}, the stable lengths of these elements are unbounded. Hence the construction can be done with the property we need.
\end{proof}
\subsection{Proof of Theorem~\ref{mthm: sparse multiplicative rigid set}}
The proof of Theorem~\ref{mthm: sparse multiplicative rigid set} is similar to the proof of Theorem~\ref{mthm: sparse rigid set}. 

Let $\Gamma$ be a group and $(X,d)$ a non-trivial hyperbolic space. Given a non-elementary or quasi-parabolic isometric action of $\Gamma$ on $X$. Let $\xi\in \Lambda_X\Gamma$. Enumerate the hyperbolic elements of $\Gamma$ as $\gamma_i$, $i\in\mathbb N^+$. Up to conjugation and changing $\gamma_i$ to $\gamma_i^{-1}$ if necessary, we may assume $<\gamma_i^-,\xi>\geq 10000i^2(\delta+1)$. Choose $\eta_i\in \Gamma$ so that $\eta_i\gamma^+\neq \gamma^-$. Let $\phi:\Gamma\to P(N)$. Define $E'_\phi$ as follows:
\[
E'_\phi=\{\gamma_i^{-\phi_{\gamma_i}(n)}\eta_i^{-1}|i,n\in\mathbb{N}^+\}.
\]

\begin{lemma}
    The set $E'_\phi$ detects the ratio of the marked length spectrum.
\end{lemma}
\begin{proof}
Let \[\overline{C}=\limsup_{\ell_d(\eta)\to \infty, \eta\in E'_\phi}\frac{\ell_d(\eta)}{\ell_{d'}(\eta)},\]
and \[\underline{C}=\liminf_{\ell_d(\eta)\to \infty, \eta\in E'_\phi}\frac{\ell_d(\eta)}{\ell_{d'}(\eta)}.\]

If $\overline{C}=\infty$ or $\underline{C}=0$, there is nothing to prove for the corresponding inequality. Hence we may assume that both of them are finite and nonzero.

  Let $(Y, d')$ be any nontrivial $\delta$-Gromov hyperbolic space with a comparable $\Gamma$ action. We identify $\partial X$ and $\partial Y$ by the homeomorphism $f$. Denote the Busemann cocycle for $X$ and $Y$ by $c_1$ and $c_2$.

Let $\gamma_i$ be a hyperbolic element. Applying Proposition~\ref{main ralation} for $\eta_i$ and $\gamma_i$, we have
\[
\limsup_{n\to\infty} |c_1(\eta_i,\gamma_i^+)+\phi_{\gamma_i}(n)\ell_d(\gamma_i)-\ell_d(\eta_i\gamma_i^{\phi_{\gamma_i}(n)})|
\leq\infty,
\]
\[
\limsup_{n\to\infty} |c_2(\eta_i,\gamma_i^+)+\phi_{\gamma_i}(n)\ell_{d'}(\gamma_i)-\ell_{d'}(\eta_i\gamma_i^{\phi_{\gamma_i}(n)})|
\leq\infty.
\]
Therefore, $\phi_{\gamma_i}(n)\ell_d(\gamma_i)=\ell_d(\eta_i\gamma_i^{\phi_{\gamma_i}(n)})+M_{n,\gamma_i}$, and $\phi_{\gamma_i}(n)\ell_{d'}(\gamma_i)=\ell_{d'}(\eta_i\gamma_i^{\phi_{\gamma_i}(n)})+M'_{n,\gamma_i}$, where both $M_{n,\gamma_i}$ and $M'_{n,\gamma_i}$ are bounded for $n$. 
Hence 
\begin{equation}\label{rr}
  \frac{\ell_d(\gamma_i)}{\ell_{d'}(\gamma_i)}=\frac{\ell_{d}(\eta_i\gamma^{\phi_{\gamma_i}(n)})+M_{n,\gamma_i}}{\ell_{d'}(\eta_i\gamma_i^{\phi_{\gamma_i}(n)})+M'_{n,\gamma_i}}  
\end{equation}

Since $\gamma_i$ is hyperbolic, $\phi_{\gamma_i}(n)\ell_d(\gamma_i)=\ell_d(\eta_i\gamma_i^{\phi_{\gamma_i}(n)})+M_{n,\gamma_i}$ implies \[\lim_{n\to \infty}\ell_d(\eta_i\gamma_i^{\phi_{\gamma_i}(n)})=\infty.\] The same is true for the $d'$ metric on $Y$. Taking limit of equation~\eqref{rr}, we have
\[
\underline{C}\leq \frac{\ell_d(\gamma)}{\ell_{d'}(\gamma)}\leq \overline{C}.
\]
This completes the proof.  
\end{proof}
\begin{proof}[Proof of Theorem~\ref{mthm: sparse multiplicative rigid set}]
    Same as in the proof of Theorem~\ref{mthm: sparse rigid set}, by taking $\phi_{\gamma_i}$ big and sparse enough, we have a sparse subset, and the only accumulation point is $\xi$. The only thing left is to show that $E'_\phi$ contains only primary elements. When the action on $X$ is proper, this follows from the same reason as in the proof of Theorem~\ref{mthm: sparse multiplicative rigid set}. If $\Gamma$ is a torsion-free hyperbolic group,  without loss of generality, we assume $\gamma$ is a primitive element in $\Gamma$. This follows from Lemma 3.6 from \cite{minasyan2010normal}. 
\end{proof}

\section{Rough cross-ratio of roughly geodesic hyperbolic space and Busemann functions}\label{Rough corss-ratio} 
In this section, we prove Theorem 5.1, which is essentially the first part of the proof of Theorem~\ref{mthm}. 
\begin{thm}\label{Thm: Gromov-product}
    Let $(X,d)$, $(Y, d')$ be two comparable nontrivial $\Gamma$-Gromov hyperbolic spaces such that $\Gamma$ is geometrically dense for $X$ and $Y$. Let $<, >_d$ and $<,>_{d'}$ be the Gromov products on $\partial X$ and $\partial Y$ for some base points, respectively.  If $\ell_d=\ell_{d'}$, then there exists $L>0$ such that
\[|<x,y>_d-<f(x),f(y)>_{d'}|\leq L\]
for all $x\neq y\in \partial X$. 
\end{thm}
The theorem is known before for cocompact actions. See Furman \cite{F}. 

\begin{proof}[Proof of Theorem~\ref{Thm: Gromov-product}]
    Since the two actions are comparable, there is a continuous homeomorphism $f:\partial X\to \partial Y$. We identify $\partial X$ and $\partial Y$ via $f$. Let $c_1$ and $c_2$ denote the Busemann cocycles of $X$ and $Y$, respectively.

    Let $\bar{c}=c_1-c_2$ and $g=<x,y>_d-<x,y>_{d'}$. The proof of Lemma~\ref{thm: marked length rigidity} shows that for fixed $\xi\in \lambda_X\Gamma=\partial X$, there exist $M>0$ such that $\bar{c}(\gamma,\xi\leq M$ for all $\gamma\in \Gamma$. By cocycle identity, $\bar{c}$ is bounded on the $\Gamma$-orbit of $\xi$. Since any $\Gamma$-orbit is dense, and the Busemann cocycles are almost continuous, we have that $\bar{c}(\gamma,\cdot)$ is uniformly bounded for $\gamma\in\Gamma$.

By Proposition~\ref{B-cocycle}, $dg$ is a bounded cocycle. Here, $dg(\gamma, (x,y))=g(x,y)-g(\gamma x,\gamma y)$. Hence, there exists a bounded function $h$ such that $dg=dh$. Indeed, we may take $h(x)=\sup_{\gamma\in \Gamma} dg(\gamma,x)$. Since $g-h$ is $\gamma$-invariant, $g-h$ is constant on any $\Gamma$-orbit. There exists a dense orbit. Because of equation~(\ref{rough continuous}), $g-h$ is bounded.  Hence, $g$ is bounded.

This completes the proof.
\end{proof}

\section{Infinite-barycenter map}
We conclude the proof of Theorem~\ref{mthm} by presenting the hyperbolic filling of a boundary. 

Consider a roughly geodesically complete Gromov hyperbolic space $(X,d)$ with a designated base point $o$. For any real numbers $K>100$ and $\delta>0$, let $\calD^R_{\delta,K}$ be the collection of pointed, $(K, \delta)$-roughly geodesically complete Gromov hyperbolic metrics $D$ on $X$, such that for all $(x,y,z,w)\in \partial^{(4)}(X)$: 
\begin{enumerate}
    \item $|<x,y>_{d}-<x,y>_{D}|<\infty.$
    \item $\sup_{y\in \partial X}\inf_{x\in \partial X}\{<x,y>_D\}\leq K-1.$
    \item $|(x,y;z,w)_d-(x,y;z,w)_D|\leq K,$
\end{enumerate}
where $(x,y;z,w)$ denotes the rough cross-ratio \eqref{rough cross-ratio}. The base point of $D$ is denoted by $o_D$, and the Gromov product for $D$ is understood with respect to the base point $o_D$.

\begin{remark}
    The space $\calD^R_{\delta, K}$ is similar to \cite{oregon2023space}*{Defnition 1.1}. $\calD^R_{\delta, K}$ is a point in $\calD(\Gamma)$ with further restrictions on the geometries since all the spaces are roughly isometric. And the metric we defined here is essentially to find the best $A$ in \cite{oregon2023space}*{Defnition 1.2}.
\end{remark}

We introduce a metric $\rho$ on $\calD^R_{\delta,K}$ defined by 
\[
\rho(D_1,D_2)=\sup_{x\neq y\in \partial X} \{|<x,y>_{D_1}-<x,y>_{D_2}|\}.
\]

For every point in $p\in X$,  there exists an element $D_p\in\calD^R_{\delta, K}$ for sufficiently large $K$. This element is the same metric space $(X,d)$ with base point $p$. It is evident from the definition and the triangle inequality that 
$\rho(D_p, D_q)\leq d(p,q)$. We denote this embedding map by $i_X$. 

Firstly, we establish that $i_X$ is a roughly geodesic embedding. Secondly, we show that the image of $i_X$ is cobounded. 

\subsection{The embedding is roughly isometric embedding}
\begin{lemma}\label{embedding}
    The embedding $i_X$ is a roughly geodesic embedding.
\end{lemma}
\begin{proof}
Consider any $p$, $q\in X$. Let $r$ be a $K$-roughly geodesic with $r(0)=p$, $r(t)=q$, $r(\infty)=x\in \partial X$ such that $|t-d(p,q)|\leq K$. Choose $y\in \partial X$ such that $<x,y>_p\geq d(p,q)$. Direct computation and equation~(\ref{sequence free}) yeild:
\[
  |(<x,y>_p-<x,y>_q)-(<y,q>_p-<x,p>_q)|\leq 8\delta.  
\]

As $r$ is roughly geodesic, applying equation~(\ref{Gromov 1; rough geodesic}) and (\ref{sequence free}) gives: 
\[<x,p>_q\leq \frac{3}{2}K+2\delta.\]
Taking a limit of the triangle inequality, we have $<x,q>_p\leq d(p,q)$. By the hyperbolic inequality (\ref{hyperolic inequality}), equation~(\ref{Gromov 2: rough geodesic}) and \eqref{rough continuous},
\[
<y,q>_p\geq \min\{<y,x>_p,<x,q>_p\}-2\delta=<x,q>_p-\delta\geq d(p,q)-\frac{5}{2}K-3\delta.
\]
Thus:
\[
\rho(D_p,D_q)\geq |<x,y>_p-<x,y>_q|\geq d(p.q)-4K-13\delta.
\]

On the other hand, for any $x\neq y\in \partial X$, consider sequences $(r_i)$ and $(s_j)$ such that: 
\[\lim_{i,j}<r_i,s_j>=<x,y>_p\] 
with $r_i\to x$, $s_j\to y$, taking the limit using this sequence, and equation (\ref{rough continuous}) along with the triangle inequality, we have 
\[
|<x,y>_p-<x,y>_q|\leq d(p,q)+4\delta.
\]
Hence
\[
\rho(D_p,D_q)\leq d(p,q)+4\delta.
\]

This completes the proof.
\end{proof}
\subsection{Coboundedness}
First, we introduce the rough version of the relative Busemann function. Given $D_1$, $D_2\in \calD^R_{\delta, K}$, let
\[
f_{D_1,D_2}(x)=\frac{1}{2}\liminf_{y\to x}<x,y>_{D_1}-<x,y>_{D_2}.
\]
Then we have $<x,y>_{d_1}-<x,y>_{d_2}-f_{D_1,D_2}(x)-f_{D_1,D_2}(y)$ is the same as
\[
\frac{1}{2}\liminf_{z\to x, w\to y}[(x,y;z,w)_{D_1}-(x,y;z,w)_{D_2}+T(z,w)]
\]
where $T(z,w)=<x,y>_{D_1}-<z,w>_{D_1}-<x,y>_{D_2}+<z,w>_{D_2}$. By equation (\ref{rough continuous}), $\liminf |T(z,w)|\leq 8\delta$.
Therefore 
\begin{equation}\label{eq: busemann}
|<x,y>_{D_1}-<x,y>_{D_2}-f_{D_1,D_2}(x)-f_{D_1,D_2}(y)|\leq \frac{1}{2}K+4\delta.
\end{equation}

We call $f_{D_1,D_2}$ the relative Busemann function of $D_1$ and $D_2$. Here are a few properties of the relative Busemann function:

\begin{lemma}\label{item 1} 
    The supremum and infimum of the Busemann function are roughly opposite.
\end{lemma}
 \begin{proof}   
    Let $x\in X$ be such that $f_{D_1,D_2}(x)\geq \sup f_{D_1,D_2}-\delta$. Then, by definition, there exist $x\neq y\in\partial X$, such that $<x,y>_{D_1}\leq K$, in view of Equation (\ref{eq: busemann}),
    \[
    -4\delta-\frac{1}{2}K+f_{D_1,D_2}(x)+f_{D_1,D_2}(y)\leq <x,y>_{D_1}-<x,y>_{D_2}\leq <x,y>_{D_1}\leq K.
    \]
    Hence $\inf(f_{D_1,D_2})\leq -\sup f_{D_1,D_2}+5\delta+\frac{3}{2}K$. A similar consideration by letting $x\neq y$ such that $f_{D_1,D_2}(y)\leq \inf f_{D_1,D_2}+\delta$, and $<x,y>_{D_2}\leq K$ leads to 
\[
    4\delta+\frac{1}{2}K+f_{D_1,D_2}(x)+f_{D_1,D_2}(y)\geq <x,y>_{D_1}-<x,y>_{D_2}\geq -<x,y>_{D_2}\geq -K.
    \]
    Hence $\inf f_{D_1,D_2}\geq -\sup f_{D_1,D_2}-\frac{3}{2}K-5\delta$. Therefore,
    \[||\inf f_{D_1,D_2}+\sup f_{D_1,D_2}||\leq 5\delta+\frac{3}{2}K.\]
\end{proof}
    
\begin{lemma}\label{item 2} 
$\rho(D_1,D_2)$ is roughly equal to $2\sup f_{D_1,D_2}$.
\end{lemma}
\begin{proof}
    
    By definition of $\rho$ and equation~(\ref{eq: busemann}), it is clear that \[\rho(D_1,D_2)\leq 2\max\{|\inf{|f_{D_1,D_2}|}|, |\sup{f_{D_1.D_2}}|\}+\frac{1}{2}K+4\delta.\] Hence, by Lemma~\ref{item 1}, 
    \begin{equation}\label{rho<2sup}
    \rho(D_1,D_2)\leq 2|\sup{f_{D_1.D_2}}|+\frac{7}{2}K+14\delta.
    \end{equation}

    On the other hand, consider a sequence $w_n\in \partial X$ converge to $x$, by (\ref{eq: busemann}),
    \[
    \begin{array}{cl}
       |f_{D_1,D_2}(x)-f_{D_1,D_2}(w_n)|\leq  & \frac{1}{2}K+4\delta+|<x,y>_{D_1}-<w_n,y>_{D_1}| \\
         & +|<x,y>_{D_2}-<w_n,y>_{D_2}|.
    \end{array}\]
    For $n$ big enough, according to equation~(\ref{rough continuous}), we have \[
    |f_{D_1,D_2}(x)-f_{D_1,D_2}(w_n)|\leq \frac{1}{2}K+8\delta.
    \]
    Hence
    $\rho(D_1,D_2)\geq |<x,w_n>_{D_1}-<x,w_n>_{D_2}|\geq 2f_{D_1,D_2}(x)-K-12\delta$ for all $x$. 
    In particular, 
    \begin{equation}\label{rho>2sup}
      \rho(D_1,D_2)\geq 2|\sup f_{D_1.D_2}|-K-12\delta.
    \end{equation}
\end{proof}

\medskip
Now we are ready to show that the image is cobounded.
\begin{lemma}
Let $(X,d)$ be a roughly $(\delta',K')$-geodesically complete Gromov hyperbolic space. Then, for $\delta$, $K$ sufficiently large, $i_X:X\to \calD^R(\delta, K)$ is a roughly isometry. Indeed,  
for any point $D\in\calD^R_{\delta,K}$, there is a point $p\in X$ so that $\rho(D,D_p)\leq 1000K+1000\delta$.
\end{lemma}
\begin{proof}
For all different points $p$, $q\in X$ and $(x,y,z,w)\in\partial^{(4)}X$. By the triangle inequality and equation~\eqref{sequence free}, taking 4 sequence converging to $x,y,z,w$ respectively, we have
\begin{enumerate}
    \item $|<x,y>_p-<x,y>_q|\leq d(p,q)+4\delta.$
    \item $|(x,y;z,w)_p-<x,y;z,w>_q|\leq 16\delta.$
\end{enumerate}
Since $X$ is $K'$-visual, there is a $K'$-geodesic connect any two point $p\in X$ and $x\in \partial X$. By the fact that $X$ is roughly complete, this $K'$ rough geodesic extends to a bi-infinite $K'$ rough geodesic. Denote the other end by $y$. By equation~\eqref{Gromov 1; rough geodesic} and \eqref{sequence free}, we have:
\[
<x,y>_p\leq \frac{3}{2}K'+2\delta.
\]
Let $\delta\geq \delta'$ and $K\geq \frac{3}{2}K'+16\delta+1$. We have $i_X:X\to \calD^R(\delta, K)$ is a rough isometric embedding by Lemma~\ref{embedding}.

Assume the image of $i_X$ is not $1000(K+\delta)$ dense, let $p\in X$ such that $1000K+1000\delta<\rho(D_p, D)<\inf\{\rho(D,D_q|q\in X)\}+\delta$. Let
\[
A=\setdef{x\in\partial X}{f_{D_p,D}\geq \frac{\rho(D_p,D)}{2}-100K-100\delta}. 
\]
By Lemma \ref{item 2}, $A$ is not empty. Taking one point $x\in A$, and from the roughly geodesic $\phi(t)$ with respected to $(X,d)$ who connected $p$ and $x$ with $\phi(0)=p$. Let $q=\phi(50K+50\delta)$. Then, for any point $y\notin A$, by definition of rough geodesics, we have
\[
f_{D_q,D}(y)\leq f_{D_p,D}(y)+\frac{1}{2}d(p,q)+\delta\leq \frac{\rho(D_p,D)}{2}-100K-100\delta+26K+26\delta.
\]
When $z\in A$, notice that $<x,z>_D\geq 0$. In view of equation~(\ref{eq: busemann}), 
\[<x,z>_{D_p}\geq \rho(D_p,D)-200K-200\delta-\frac{1}{2}K-4\delta.\]

By the definition of the relative Busemann function, we have:
\[
f_{D_q,D}(z)\leq f_{D_p,D}(z)+\frac{1}{2}\limsup_{w\to z}<z,w>_q-<z,w>_p.
\]
It is a direct computation and equation~(\ref{sequence free}) that shows that:
\[<z,w>_q-<z,w>_p\leq<w,p>_q-<z,q>_p+8\delta.\] 

First, since $x$ is the limit of a rough geodesic, equation~(\ref{sequence free}) and
\[K+2\delta\geq(x,p)_q\geq\min\{<w,p>_q,<x,w>_q\}-2\delta,\]
implies that $<w,p>_q\leq K+4\delta$. Second, 
\[<z,q>_p\geq \min\{(x,z)_p,<x,q>_p\}-2\delta.\]
Therefore, $<z,q>_p\geq d(p,q)-K-2\delta=48K+48\delta$ from equation~(\ref{Gromov 1; rough geodesic}).
It follows that
\[
f_{D_q,D}(z)\leq f_{D_p,D}(z)-\frac{47}{2}K-18\delta.
\]
Therefore by equation~(\ref{rho>2sup}), 
\[
\sup_{z\in A} F_{D_q,D}(z)\leq \frac{\rho(D_p,D)}{2}+\frac{1}{2}K+6\delta-\frac{47}{2}K-18\delta.
\]
And recall for $y\notin A$,
\[
\sup_{y\notin A}F_{D_q,D}(y)\leq \frac{\rho(D_p,D)}{2}-74K-74\delta.
\]
Finally, according to equation~(\ref{rho<2sup}),
\[
\rho(D_q,D)\leq 2\sup F_{D_q,D}+3.5K+14\delta\leq\rho(D_p,D)-42.5K-10\delta.
\]
Contradict to the choice of $p$. Hence, the image of $X$ is $1000K+1000\delta$ dense.

This completes the proof.
\end{proof}
\subsection{Proof of Theorem~\ref{mthm}}
\begin{proof}
Choose $\delta$, $K$ bigger enough so that $i_X$ and $i_Y$ are rough isometrics. Since both maps are $\Gamma$-equivalent. The result follows. 
\end{proof}
\section{Proof of Theorem~\ref{Mthm: coset}}
\begin{proof}

Since $\ell_d(h\gamma)=\ell_{d}(\gamma^{-1}h^{-1})$. Up to switching $h$, we may assume the marked length spectrum agrees on a left coset $h\Gamma'$. 

Same as in the proof of Theorem~\ref{mthm}, we identify $\partial X$ and $\partial Y$ via the comparable map $f$. Let $c_1$ and $c_2$ denote the Busemann cocycle of $X$ and $Y$, respectively. Let $c=c_1-c_2$ and $g(x,y)=<x,y>_d-<x,y>_{d'}$. Then form section~\ref{Busemann cocycle}, we have:
    \begin{equation}\label{eq:11}
        c(h\gamma,(h\gamma)^{\pm})=0,
    \end{equation} 
    \begin{equation}\label{eq:12}
        |g(x, y)-g(\eta x,\eta y)-c(\eta,x)-c(\eta,y)|\leq 8\delta.
    \end{equation} 
for all $\gamma\in \Gamma'$, $\eta\in\Gamma$, $x$, $y\in\partial X$.

    First, we show that $\ell_d=\ell_{d'}$ on $\Gamma'$. By Proposition~\ref{Prop: stable length}, it is equivalent to that for all hyperbolic element $\gamma\in \Gamma'$,
    \[
    \lim_n\frac{c(\gamma^n,\gamma^-)}{n}=0.
    \]
    
    Let $\gamma\in\Gamma'$ be a hyperbolic element so that $h\gamma^+\neq \gamma^-$. Then by equation~$\eqref{eq:12}$,
    \begin{equation}\label{main}
        |g(\gamma^-,(h\gamma^n)^+)-g(h\gamma^{-},(h\gamma^n)^+)-c(h\gamma^n, \gamma^-)-c(h\gamma^n,(h\gamma^n)^+)|\leq 8\delta.
    \end{equation}
    
    Now we study each term in equation~\eqref{main}, and connect it to $c(\gamma^n,\gamma^-)$.
 \begin{enumerate}
     \item In view of Lemma~\ref{S} and equation~\eqref{rough continuous}, we have
     \[
     \limsup_n|g(\gamma^-,(h\gamma^n)^+)|\leq |g(\gamma^-,h\gamma^+)|+4\delta.
     \]
     \item For the same reason,
     \[
     \limsup_n g(h\gamma^-,(h\gamma^n)^+)\leq |g(h\gamma^-,h\gamma^+)|+4\delta.
     \]
     \item By roughly cocycle identity (\ref{almost cocycle}) 
     \[
     |c(h\gamma^n,\gamma^-)-c(h,\gamma^-)-c(\gamma^n,\gamma^-)|\leq 4\delta.
     \]
     Since by the triangle inequality, $c(h,\gamma^-)\leq d(o,ho)+d'(o',ho')$ where $o$ and $o'$ are base point of $X$ and $Y$, respectively. Hence, $c(h\gamma^n,\gamma^-)$ is uniformly close to $c(\gamma^n,\gamma^-)$.
     \item By equation~\ref{eq:11}, 
     \[
     c(h\gamma,(h\gamma)^{\pm})=0.
     \]
 \end{enumerate}
    Put them back in equation~\eqref{main}, we have, for $n$ big enough, 
    \[
    |c(\gamma^n,\gamma^-)|\leq 20\delta+d(o,ho)+d'(o',ho')+|g(\gamma^-,h\gamma^+)|+|g(h\gamma^-,h\gamma^+)|+1.
    \]
    Hence $\lim_n \frac{c(\gamma^n,\gamma^-)}{n}=0$.

    Now, view $X$ and $Y$ as a $\Gamma'$-space, by proof of Theorem~\ref{mthm}. Both $i_X:X\to \mathcal{D}_{\delta,K}$ and $i_Y:Y \to \mathcal{D}_{\delta,K}$ are rough isometics for some sufficient big $K$. Since both $\Gamma$ actions preserve the rough cross-ratio, $i_X$ and $i_Y$ are (roughly) $\Gamma$-equivariant. Hence, the result follows.  
\end{proof}

\end{document}